\documentclass[11pt]{amsart}
\usepackage{amssymb}
\usepackage{url}
\usepackage{enumerate}
\usepackage{tikz}
\newtheorem{lemma}{\bf Lemma}
\newtheorem{theorem}[lemma]{\bf Theorem}

\newtheorem{proposition}[lemma]{\bf Proposition}

\DeclareMathOperator{\id}{id}

\begin{document}
\parskip = 0mm
\title[A Problem of J\'onsson and McKenzie from 1982]{Another Problem of J\'onsson and McKenzie from 1982\\ \small Refinement Properties for Connected Powers of Posets}
\author{Jonathan David Farley}
\address{Department of Mathematics, Morgan State University, 1700 E. Cold Spring Lane, Baltimore, MD 21251, United States of America, {\tt lattice.theory@gmail.com}}



\begin{abstract}
In 1982, J\'onsson and McKenzie posed the following problem: ``Find counter examples (or prove that none exist) to the refinement of $A^C\cong B^D$ [$A$, $B$, $C$, and $D$ non-empty posets] under'' the condition ``$C$, $D$, and $A^C$ are finite and connected.'' That is, in this situation, are there posets $E$, $X$, $Y$, and $Z$ such that $A\cong E^X$, $B\cong E^Y$, $C\cong Y\times Z$, and $D\cong X\times Z$?

In this note, this problem is solved.
\end{abstract}

\thanks {The author would like to thank Dr. Bernd S. W. Schr\"oder for suggestions on improving this note.}

\maketitle


\def\Qa{\mathbb{Q}_0}
\def\Qb{\mathbb{Q}_1}
\def\Q{\mathbb{Q}}
\def\card{{\rm card}}
\parskip = 2mm
\parindent = 10mm
\def\Part{{\rm Part}}
\def\P{{\mathcal P}}
\def\Eq{{\rm Eq}}
\def\cld{Cl_\tau(\Delta)}
\def\Csing{{\mathcal C}_{\{*\}}}
\def\Cftwo{{\mathcal C}_{{\rm fin}\rangle1}}
\def\Cinf{{\mathcal C}_{\infty}}
\def\Pcf{{\mathcal P}_{\rm cf}}
\def\Fn{{\mathcal F}_n}
\def\proof{{\it Proof. }}


\vspace*{-4mm} 

Let $E$ and $X$ be posets.  Define $E^X$ to be the poset of order-preserving maps from $X$ to $E$ where, for $f,g\in E^X$, $f\le_{E^X} g$ if $f(x)\le_E g(x)$ for all $x\in X$.

\begin{center}

    \begin{tikzpicture}[scale=.85]

    \draw[fill] (-5,0) circle (.05cm);
    \draw (-5,0) -- (-5,1);
    \draw[fill] (-5,1) circle (.05cm);
    \draw (-5,0) -- (-4,1);
    \draw[fill] (-4,1) circle (.05cm);

    \draw[fill] (-4,0) circle (.05cm);
    \draw (-4,0) -- (-4,1);
    \draw[fill] (-4,1) circle (.05cm);
    \draw (-4,0) -- (-5,1);

    \draw[fill] (-3.5,1) circle (.05cm);
    \draw (-3.5,1) -- (-3.5,2);
    \draw[fill] (-3.5,2) circle (.05cm);

    \draw (-2,1) node {\large $\cong$};

    \draw[fill] (-1,0) circle (.05cm);
    \draw (-1,0) -- (-1,1);
    \draw[fill] (-1,1) circle (.05cm);
    \draw (-1,1) -- (-1,2);
    \draw[fill] (-1,2) circle (.05cm);
    \draw (-1,0) -- (0,1);
    \draw[fill] (0,1) circle (.05cm);

    \draw[fill] (2,0) circle (.05cm);
    \draw (2,0) -- (2,1);
    \draw[fill] (2,1) circle (.05cm);
    \draw (2,1) -- (2,2);
    \draw[fill] (2,2) circle (.05cm);
    \draw (2,0) -- (1,1);
    \draw[fill] (1,1) circle (.05cm);

    \draw (0,1) -- (2,2);
    \draw (1,1) -- (-1,2);

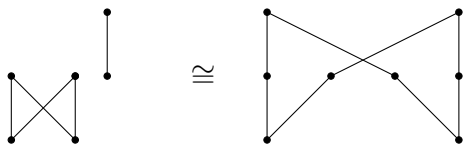
    \draw (-1,-1) node {\bf Figure. \sl $E^X$ where $E$ is the 4-element crown and $X$ the 2-element chain.};

    \end{tikzpicture}
    
\end{center}

\noindent Many of the usual laws of arithmetic hold for the arithmetic of ordered sets, with one unusual one: if $A$ is connected and non-empty, then $(U+S)^A\cong U^A+S^A$ \cite[Proposition 4.1]{McKJC}.  So if $E$, $X$, $Y$, and $Z$ are posets,
$$
(E^X)^{Y\times Z}\cong E^{X\times Y\times Z}\cong(E^Y)^{X\times Z}\quad\text{($\natural$)}
$$
What if we have posets $A$, $B$, $C$, and $D$ such that $A^C\cong B^D$?  Can we find ``refining'' posets $E$, $X$, $Y$, and $Z$ such that the isomorphism is naturally explained (by $\natural$), where
$$
A\cong E^X\text{,}\ B\cong E^Y\text{,}\ C\cong Y\times Z\text{, and}\ D\cong X\times Z\text{?}
$$

In 1982, Bjarni J\'onsson (an invited speaker at the 1974 International Congress of Mathematicians \cite[\S3]{NatAH}) and Ralph McKenzie, a professor at the University of California at Berkeley, posed the following problem \cite{JonMcKHB}:

``{\large P}ROBLEM 12.1. Find counter examples (or prove that none exist) to the refinement of $A^C\cong B^D$ under any of the conditions:\dots

\noindent (ii) $C$, $D$, and $A^C$ are finite and connected.''

We solve this problem.\footnote{Since J\'onsson and McKenzie used the word ``any'' (they listed five sets of conditions), we assert that, technically, we have solved their Problem 12.1.}

For posets $E$ and $X$, let $\mathcal D(E^X)$ denote
$$
\{g\in E^X\mid g\ \text{is constant on each connected component of}\ X\}
$$
and for $X\ne\emptyset$ let $\mathcal C(E^X)$ denote
$$
\{f\in E^X\mid f\ \text{is in the same connected component as}\  g\ \text{for some}\ g\in\mathcal D(E^X)\}
$$
\noindent Recall that a subset $X$ of a poset $E$ is a {\it retract} if there exists an order-preserving map (a {\it retraction}) $\rho:E\to X$ such that $\rho\restriction_X=\id_X$.  A general reference is \cite{SchJC}.

\begin{lemma} Let $E$ and $X$ be posets such that $X\ne\emptyset$.  Then $E$ is order-isomorphic to a retract of $\mathcal C(E^X)$.
\end{lemma}

\proof For all $e\in E$, the constant map from $X$ to $E$ with image $\{e\}$, $\langle e\rangle$, is in $\mathcal D(E^X)\subseteq\mathcal C(E^X)$.  Fix $x_0\in X$.  Define a map
$$
\rho:\mathcal C(E^X)\to\{\langle e\rangle\mid e\in E\}
$$
\noindent by $\rho(f)=\langle f(x_0)\rangle$ for all $f\in\mathcal C(E^X)$.

If $f,g\in\mathcal C(E^X)$ and $f\le g$, then $f(x_0)\le g(x_0)$, so $\langle f(x_0)\rangle\le\langle g(x_0)\rangle$.  Also, if $e\in E$, then $\rho(\langle e\rangle)=\big\langle\langle e\rangle(x_0)\big\rangle=\langle e\rangle$. $\qed$

\begin{lemma} Let $E$ and $X$ be posets such that $X\ne\emptyset$.  If $E^X$ is connected, then $\mathcal C(E^X)=E^X$. $\qed$
\end{lemma}

\begin{lemma} Let $A$ and $S$ be posets.  Let $\rho: A\to I$ be a retraction onto a subset $I$ of $A$.  Define $\sigma: A^S\to I^S$ for all $f\in A^S$ by $\sigma(f)=\rho\circ f$.  Then $\sigma$ is a retraction.
\end{lemma}

\proof Let $f,g\in A^S$ be such that $f\le g$.  Then for all $s\in S$, $f(s)\le g(s)$, so
$$
[\sigma(f)](s)=(\rho\circ f)(s)=\rho\big(f(s)\big)\le \rho\big(g(s)\big)=(\rho\circ g)(s)= [\sigma(g)](s).
$$
\noindent Hence $\sigma(f)\le\sigma(g)$.

If $f\in I^S$, then for all $s\in S$,
$$
[\sigma(f)](s)=(\rho\circ f)(s)=\rho\big(f(s)\big)=f(s)\text{,}
$$
\noindent so $\sigma(f)=f$. $\qed$

\begin{proposition} Let $S$, $P$, $Q$, and $R$ be posets such that $Q,R,S\ne\emptyset$, and $P\cong\mathcal C(Q^R)$.  Then $Q^S$ is order-isomorphic to a retract of $P^S$.  Hence, if $P^S$ is connected, so is $Q^S$.
\end{proposition}

\proof Apply Lemma 3 to the retraction of Lemma 1.

The image of a connected poset under an order-preserving map is connected. $\qed$

\begin{lemma} Let $A$, $S$, and $U$ be non-empty posets such that $U^{S\times A}$ is connected.  Then $U^S$ is connected.
\end{lemma}

\proof Use \cite[Proposition 4.1(10)]{McKJC} and the fact that $U^{S\times A}\cong (U^S)^A$. $\qed$

\begin{theorem} Let $A$, $B$, $C$, and $D$ be non-empty posets such that $C$, $D$, and $A^C$ are finite and connected and $A^C\cong B^D$.  Then there exist posets $E$, $X$, $Y$, and $Z$ such that $A\cong E^X$, $B\cong E^Y$, $C\cong Y\times Z$, and $D\cong X\times Z$.
\end{theorem}

\proof Obviously $A$ and $B$ must be finite. Since $C$ and $D$ are non-empty and connected and $A^C$ is connected, then $A$, $B$, and $B^D$ are connected \cite[Proposition 4.1(10)]{McKJC}.  By \cite[Theorem 4]{FarBJc} and Lemma 2, we know that there exist finite, non-empty, connected posets $E$, $X$, $Y$, and $Z$ such that $A\cong\mathcal C(E^X)$, $B\cong\mathcal C(E^Y)$, $C\cong Y\times Z$, and $D\cong X\times Z$---i.e., $\mathcal C(E^X)^{Y\times Z}\cong A^C\cong B^D\cong\mathcal C(E^Y)^{X\times Z}$.  By Lemma 5, $\mathcal C(E^X)^Y$ and $\mathcal C(E^Y)^X$ are connected, so, by Proposition 4, $E^Y$ and $E^X$ are connected.  Thus, by Lemma 2, $A\cong\mathcal C(E^X)=E^X\ \text{and}\ B\cong\mathcal C(E^Y)=E^Y$. $\qed$
 
It would be interesting to find a ``strict'' version of Theorem 6 (cf. \cite[Proposition 3.1]{McKJC}).

\medskip

In \cite{KolHH}, Kolibiar defines congruences on connected posets and proves they form a lattice.  He told the author that he did not know if his congruence lattices were distributive \cite{KolI?}.  But if they are, perhaps that can be used to prove the refinement theorem (cf. \cite[Chapter 5]{McKMcNTayHG}).


\end{document}